# ADAPTIVE MULTISCALE DETECTION OF FILAMENTARY STRUCTURES IN A BACKGROUND OF UNIFORM RANDOM POINTS[1]

By Ery Arias-Castro, David L. Donoho and Xiaoming Huo

*University of California, San Diego, Stanford University and Georgia Institute of Technology*

We are given a set of $n$ points that might be uniformly distributed in the unit square $[0, 1]^2$. We wish to test whether the set, although mostly consisting of uniformly scattered points, also contains a small fraction of points sampled from some (a priori unknown) curve with $C^\alpha$-norm bounded by $\beta$. An asymptotic detection threshold exists in this problem; for a constant $T_-(\alpha, \beta) > 0$, if the number of points sampled from the curve is smaller than $T_-(\alpha, \beta)n^{1/(1+\alpha)}$, reliable detection is not possible for large $n$. We describe a multiscale significant-runs algorithm that can reliably detect concentration of data near a smooth curve, without knowing the smoothness information $\alpha$ or $\beta$ in advance, provided that the number of points on the curve exceeds $T_*(\alpha, \beta)n^{1/(1+\alpha)}$. This algorithm therefore has an optimal detection threshold, up to a factor $T_*/T_-$.

At the heart of our approach is an analysis of the data by counting membership in multiscale multianisotropic strips. The strips will have area $2/n$ and exhibit a variety of lengths, orientations and anisotropies. The strips are partitioned into anisotropy classes; each class is organized as a directed graph whose vertices all are strips of the same anisotropy and whose edges link such strips to their "good continuations." The point-cloud data are reduced to counts that measure membership in strips. Each anisotropy graph is reduced to a subgraph that consist of strips with significant counts. The algorithm rejects $\mathbf{H}_0$ whenever some such subgraph contains a path that connects many consecutive significant counts.

## 1. Introduction.

Received July 2003; revised March 2005.

[1]Supported in part by NSF Grants DMS-00-77261, DMS-01-40587 (FRG) and DMS-05-05303, ONR MURI, and a contract from DARPA ACMP.

*We cannot help but see faces and castles in clouds, monsters in ink-blots and exotic forms in random dots. Form is so central to human perception that, I am told, it is extremely difficult to prove something random or formless. Mae-Wan Ho* [20].

Suppose we have $n$ data points $X_i \in [0,1]^2$ which at first glance seem uniformly distributed in the unit square. On cursory visual inspection, it seems that a suspiciously large number of the data points fall along a smooth curve. However, the curve on which these points lie has only been identified *after* inspection of the data. We know that the human visual system has the ability to "hallucinate" curvilinear structure in truly random point clouds. We are therefore concerned about the reliability of the perceived pattern and wish to follow an objective procedure for testing the existence of filamentary structure: this procedure should reliably separate filamentary structure from random scatter and be computationally tractable.

This is a prototype for various practical imaging problems that range from surveillance to road and streamed tracking to particle physics [1, 31, 34]. In all cases, the observer is looking for evidence of a filamentary structure in a background of heavy clutter.

As a first attempt to formalize matters, consider the problem of testing

$$\mathbf{H}_0 : X_i \overset{\text{i.i.d.}}{\sim} \text{Uniform}(0,1)^2,$$

versus

$$\mathbf{H}_1(\alpha, \beta) : X_i \overset{\text{i.i.d.}}{\sim} (1 - \varepsilon_n) \text{Uniform}(0,1)^2 + \varepsilon_n \text{Uniform}(\text{graph}(f)),$$

where $f \in \text{Hölder}(\alpha, \beta)$ is unknown. Here, for $1 < \alpha \leq 2$, Hölder$(\alpha, \beta)$ is the class of functions $g : [0,1] \to [0,1]$ with continuous derivative $g'$ that obeys $|g'(x) - g'(y)| \leq \alpha \beta |x - y|^{\alpha - 1}$. In words, we believe that a relatively small fraction $\varepsilon_n$ of points lie on a smooth curve in the plane.

1.1. *"Connect the dots"*. In our previous work [5], it was shown that when $\alpha$ and $\beta$ are fixed and known, there is a detector based on the principle that, under $\mathbf{H}_0$, no Hölder$(\alpha, \beta)$ curve can pass through a very large number of points in a random point cloud. More particularly, we know that there is a threshold $T_+ = T_+(\alpha, \beta)$ such that:

- If $T < T_+$, we have that, with probability tending to 1, *there exists* a Hölder$(\alpha, \beta)$ curve that contains at least $T \cdot n^{1/(1+\alpha)}$ points (out of $n$).
- If $T > T_+$, we have that, with probability tending to 1, there *does not exist* a Hölder$(\alpha, \beta)$ curve that contains more than $T \cdot n^{1/(1+\alpha)}$ points $(X_i)_{i=1}^n$.

(More concretely, if we deal with Lipschitz curves with |slope| $\leq 1$, we have found empirically that for moderate $n \approx 1000$, there will frequently be some Lipschitz curve that contains $\sqrt{n}$ data points, but rarely will there be one

that contains more than $3\sqrt{n}$ points.) So, if we happen to notice a curve passing through substantially more than $T_+ \cdot n^{1/(1+\alpha)}$ points, we have a strong basis to reject the null hypothesis of pure randomness. Moreover, to within a constant factor this threshold is optimal; no sequence of tests can be reliable for detecting substantially fewer than $T_- \cdot n^{1/(1+\alpha)}$ points, for a certain $T_- > 0$.

Elaborating this connect-the-dots (CTD) principle leads to a formal hypothesis test based on searching for curves that contain large numbers of points. Let $N_n(A) = \#\{\{X_i\} \cap A\}$ denote the measure that counts how many points lie in the set $A$. Searching for a curve that contains the maximal number of points leads to the optimization problem

$$N_n^*(\alpha, \beta) = \max\{N_n(\mathrm{graph}(f)) : f \in \mathrm{H\ddot{o}lder}(\alpha, \beta)\},$$

which searches over all Hölder$(\alpha, \beta)$ graphs and rejects $\mathbf{H}_0$ for values of $N^*$ that substantially exceed $T_+ \cdot n^{1/(1+\alpha)}$.

The CTD approach, while very instructive, does not address the concerns of someone who might actually be interested in performing a test on real data. Such concerns include:

- *Computational burden.* The task of finding the largest number of points on a Hölder$(\alpha, \beta)$ curve seems to us to be computationally impractical unless $\alpha \in \{1, 2\}$.
- *Unknown $\alpha$, $\beta$.* The CTD approach assumes a specific $\alpha$-$\beta$ combination. Instead, we desire an algorithm that works regardless of the specific values of $\alpha \in (1, 2]$ and $\beta > 0$.
- *Fragments.* The CTD approach searches only for graphs that extend all the way across the square from $x = 0$ to $x = 1$. Instead, one wants an algorithm that works even for short graphs.
- *General planar curve.* The CTD approach assumes that the underlying curve can be parametrized as a graph. It seems important to search for general curves rather than just graphs—for example, curves that loop around in the plane.

1.2. *An adaptive multiscale approach.* In this paper we describe an approach that addresses the concerns just listed. Our proposal:

1. Works across for a wide range of $(\alpha, \beta)$, and only requires knowing a bound on the maximum slope of the curve.
2. Detects the presence of $\mathbf{H}_1$ provided

$$\varepsilon_n > T_* \cdot n^{-\alpha/(\alpha+1)}$$

for a constant $T_*$ which depends on $\alpha$, $\beta$ and other factors. In view of earlier results, this is optimally sensitive to within a factor $T_*/T_-$.

3. Runs in $O(n^2 \cdot \log(n))$ flops.
4. Extends naturally to detect general planar curves that are not graphs.
5. Extends naturally to detect target filaments of unknown extent that, in
   large samples, can be very short compared to the image extent.

The detector is based on a kind of *multiscale geometric analysis* of the
data set, using a multiscale dictionary of parallelogram strips that exhibit a
variety of lengths, locations, orientations and aspect ratios. The idea is to
count membership of data points in various strips, to identify strips with
significantly large counts and to search for long runs of significantly large
counts in collections of strips that are "good continuations" of each other.

The detector is adaptive to the unknown smoothness $(\alpha, \beta)$ in the sense
that it achieves near-optimal performance over a wide range: $1 < \alpha \le 2$,
$\beta > 0$. (This notion of adaptivity parallels the notion of adaptive near-
minimaxity in nonparametric smoothing, in which a single estimator, able to
perform in a near-minimax way across a whole range of different smoothness
conditions, is called *adaptive to unknown smoothness* [17].) Ultimately, such
adaptivity flows from ideas behind Lemma 2.2 below, which show that our
class of strips has certain covering properties uniformly over each smooth-
ness class in the range $1 < \alpha \le 2$.

An interesting aspect of our approach is how simply and naturally the
principle of good continuation appears and leads to a solution.

Note that here we consider only the *presence* of the underlying curve—a
*detection* problem. Another question—the *estimation* problem—is to locate
the position of the curve accurately. The performance of our procedure for
estimation will *not* be addressed here.

1.3. *Contents.* Section 2 describes our underlying multiscale data struc-
tures. Section 3 describes our adaptive algorithm in general terms and gives
a statement of our main result. Section 4 describes the threshold settings
that underlie our algorithm, while Sections 5 and 6 analyze its behavior un-
der $\mathbf{H}_0$ and $\mathbf{H}_1$, respectively. Section 7 finishes the paper with a discussion
of related work.

**2. Multiscale anisotropic strips and good continuation.** Our data struc-
tures comprise a multiscale collection of anisotropic, tilted planar regions
and a sequence of directed graphs that organize them. We use ideas and
notation common in dyadic multiscale analysis (e.g., dyadic partitioning)
[10, 14, 15, 16, 17, 24]; in particular, we assume that $n$ is large and find it
convenient to let $J = \lceil \log_2(n) \rceil$ denote its dyadic logarithm. The variable $j$
will index dyadic scales $2^{-j}$ and will range through $0 \le j \le J$.

In our construction we fix in advance $S > 1$ (e.g., 2 or 4); this controls
the maximum |slope| we will be able to detect.

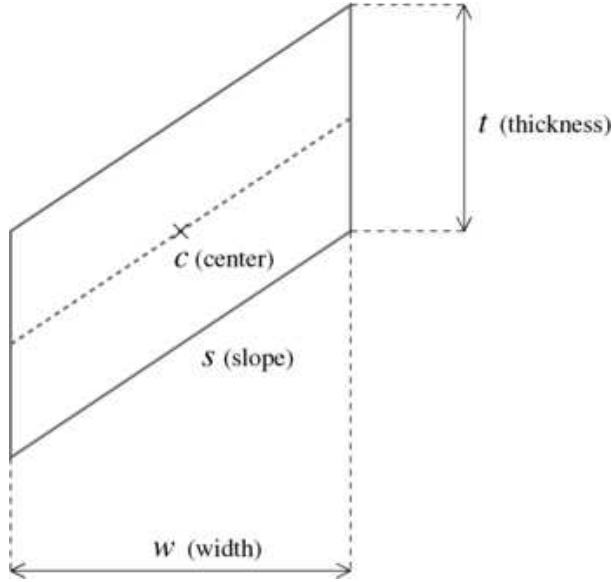

Fig. 1. *An anisotropic strip $R$.*

Let $R(j, k, \ell_1, \ell_2)$ be a parallelepiped with vertical sides that is $w = 2^{-j}$ wide by $t = 2^{-(J-j)+1}$ thick. Here $j$ runs through our set of scale indices $\{0, \ldots, J\}$. For examples, see Figure 1. The regions in question have a midline that bisects them vertically and will be tilted (sheared) at a variety of angles. Notice that these regions are highly anisotropic. While the whole collection implicitly depends on $n$, we suppress this in our notation. Moreover, the width $w$ and thickness $t$ depend on $j$ and $n$, but we also suppress this in our notation. Note that the degree of anisotropy is the same for all regions that share a common value of $j$; we generally focus only on one anisotropy class $j$ at a time.

The parameters $k$ and $\ell_i, i = 1, 2$, control the horizontal location of the regions and the vertical location and slope of the midline. There is an underlying assumption that we are interested only in regions whose major axis has a slope bounded in absolute value by $S$. The mapping between these discrete parameters is intended to insure that the regions pack together horizontally and that they are fairly closely spaced in both vertical position and slope. Let $\delta_1 = t/4$ and $\delta_2 = t/(4w)$ (these again depend implicitly on $j$ and $n$). The parallelepiped $R(j, k, \ell_1, \ell_2)$ will be centered at $c = ((k + 1/2)w, \ell_1 \delta_1)$ and its midline will have slope $s = \ell_2 \delta_2$. Here $0 \le k < w^{-1}$, $\ell_1$ runs through the set $0, \ldots, \delta_1^{-1} - 1$ and $\ell_2$ runs through the set $-S\delta_2^{-1}, \ldots, S\delta_2^{-1}$.

We gather all such regions at level (scale) $j$ in $\mathcal{R}(j) = \{R(j, k, \ell_1, \ell_2) : k, \ell_1, \ell_2\}$.

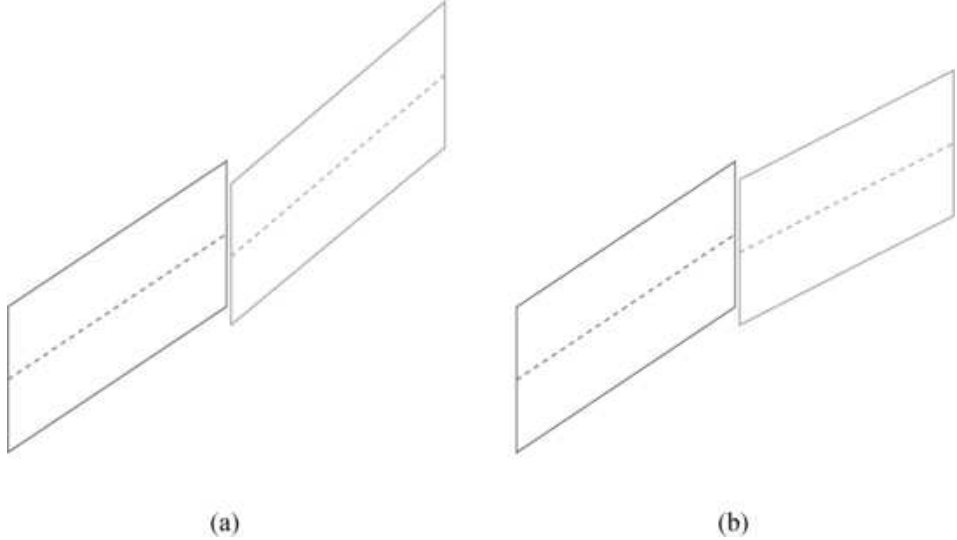

Fig. 2. *Examples of* good *continuations; the midlines have about the same slope and there is a substantial part of a side in common.*

To organize the regions, we define a directed graph $\mathcal{G}(j) = (\mathcal{V}(j), \mathcal{E}(j))$, with vertices $\mathcal{V}(j)$ and edges $\mathcal{E}(j)$. The vertices are simply the regions in $\mathcal{R}(j) \colon \mathcal{V}(j) \equiv \mathcal{R}(j)$. The edges connect regions to their *good continuations*, namely, regions that are horizontally adjacent, and that have altitudes and slopes that are nearly the same—less than $\delta_1$ and $\delta_2$ apart, respectively. Formally, we have the directed edges in $\mathcal{E}(j)$,

$$(k, \ell_1, \ell_2) \mapsto (k+1, \ell_1 + \ell_2 + u, \ell_2 + v),$$

where $|u| \le 4$, $|v| \le 4$.

Figures 2 and 3 illustrate good continuation and bad continuation, respectively.

This graphical structure, while very simple, has a perhaps surprising property: it allows us efficiently to cover the graphs of general smooth functions that exhibit any of a range of smoothnesses. This claim is summarized in three lemmas. The first lemma associates to each Hölder class a specific anisotropy graph.

LEMMA 2.1. *For each fixed $(\alpha, \beta)$ combination with $1 \le \alpha \le 2$, we have that for all sufficiently large $n$ there is $j^* = j^*(\alpha, \beta; n)$ so that $w = w(j^*)$ and $t = t(j^*)$ obey*

$$(2.1) \qquad\qquad 2\beta w^\alpha \le t < 16\beta w^\alpha.$$

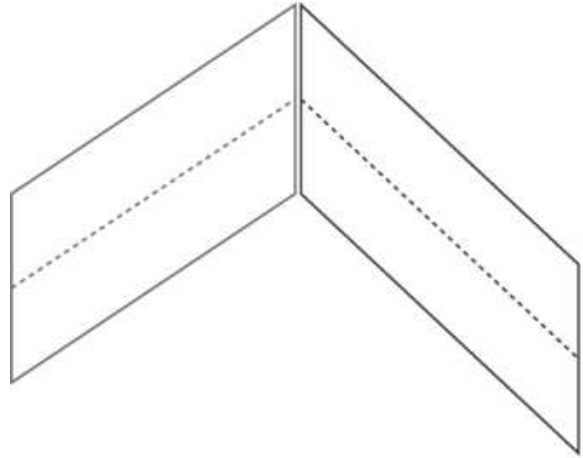

**(a)**

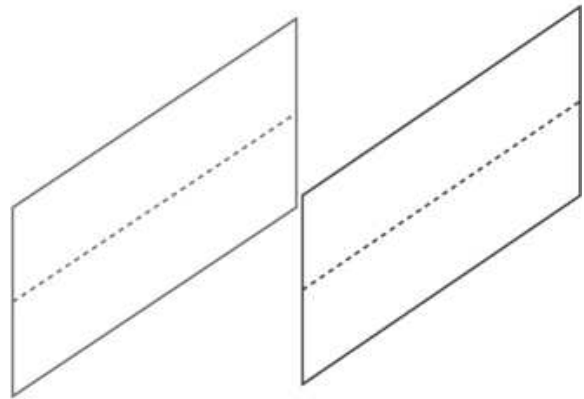

**(b)**

FIG. 3.   *Examples of* bad *continuations; either the midlines have very different slopes or the sides are effectively disjoint.*

The next result shows that regions in the anisotropy class $\mathcal{R}(j^*(\alpha,\beta))$ are well *adapted* to cover fragments of the graphs of the associated Hölder$(\alpha,\beta)$ class.

LEMMA 2.2.   *Let* $j = j^*(\alpha,\beta)$ *and suppose* $f$ *is a Hölder$(\alpha,\beta)$ function with a domain that contains* $I_k = [kw, (k+1)w)$. *Set* $x_k = (k+1/2)w$, *let* $\ell_{1,k}\delta_1$ *be the closest multiple of* $\delta_1$ *to* $f(x_k)$ *and let* $\ell_{2,k}\delta_2$ *be the closest multiple of* $\delta_2$ *to* $f'(x_k)$. *We say that the region* $R(j, k, \ell_{1,k}, \ell_{2,k})$ *is associated to* $f$ *on* $I_k$. *This strip covers the graph of* $f$ *over their common domain:*

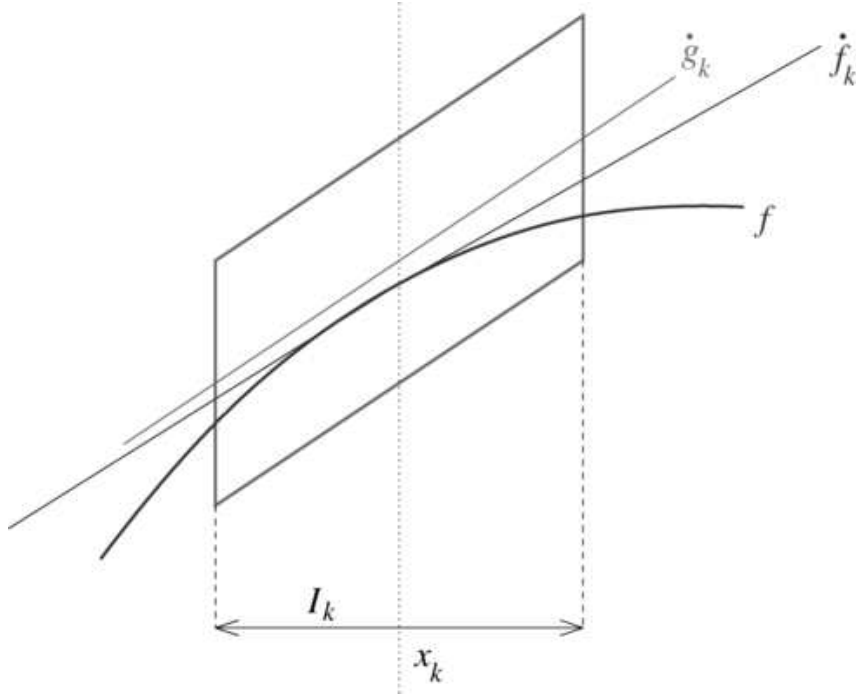

Fig. 4. *Graph of $f$ covered by its $k$th associated region $R_k$ at scale $j$; $\dot{f}_k$ is the tangent to* graph$(f)$ *at* $x_k$ *and* $\dot{g}_k$ *is the midline of* $R_k$.

$$\text{graph}(f|I_k) \subset R(j, k, \ell_{1,k}, \ell_{2,k}). \tag{2.2}$$

(See Figure 4.)

The final lemma in the sequence shows that every function in the Hölder$(\alpha, \beta)$ class corresponds to a covering sequence of regions that makes a connected path in $\mathcal{G}(j)$.

LEMMA 2.3. *Let* $j = j^*(\alpha, \beta)$ *and suppose* $f$ *is a Hölder$(\alpha, \beta)$ function on* $[0, 1]$. *For each* $k = 0, \ldots, w^{-1} - 1$, *consider the region* $R_k \equiv R(j, k, \ell_{1,k}, \ell_{2,k})$ *associated to* $f$ *by the procedure mentioned in Lemma* 2.2. *The sequence of strips* $\mathcal{T}_j(f) \equiv \{R_k : 0 \leq k < w^{-1}\}$ *consists of spatially adjacent regions, making a kind of tube. When viewed as vertices of* $\mathcal{G}(j)$, *the* $(R_k)$ *are neighbors in* $\mathcal{G}(j)$, *that is,* $R_k$ *and* $R_{k+1}$ *can be connected using edges in* $\mathcal{E}(j)$. *Therefore,* $\mathcal{T}_j(f)$ *corresponds to a path in* $\mathcal{G}(j)$.

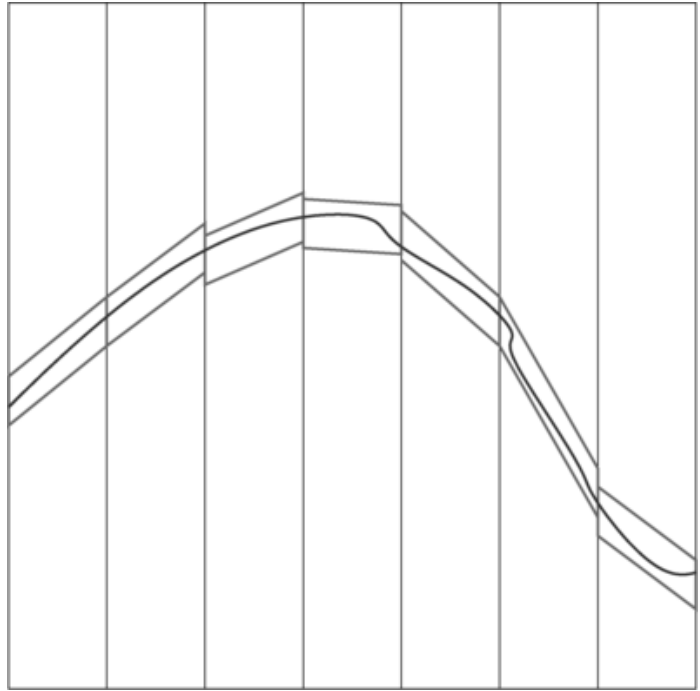

Fig. 5. *Graph of $f$ covered by its associated tube $\mathcal{T}_j(f)$ at scale $j$.*

These lemmas together show that, while the graphical structure itself is based on very simple rules, it is able to associate paths in the graph with custom-fitting tubes that cover the graphs of very different kinds of smooth functions. The proofs of these lemmas are given in the Appendix. Figure 5 illustrates the idea.

**3. The multiscale significant-runs algorithm.** We now describe the complete algorithm for analysis of point-cloud data $(X_i)$ looking for suspected curvilinear structure. It depends on a counting threshold $N^*$ and a length threshold $L_n^*$, both to be defined later. The algorithm has several steps:

1. *Counting membership in anisotropic strips.* For every region $R$, in every anisotropy class, we count the number of data that fall into that region,

$$N(R) = \#\{i : X_i \in R\}.$$

2. *Identifying significant counts.* We define a *significance indicator*, which is nonzero when the counts exceed a threshold,

$$s(R) = \mathbb{1}_{\{N(R) > N^*\}}.$$

The significance indicator may be viewed as a label on the regions $R$, producing a sequence of a labeled graphs

$$\Sigma(j) = (\mathcal{V}(j), \mathcal{E}(j), \sigma(j)),$$

where $\sigma(j) = (s(R))$ gives the labels on $R \in \mathcal{R}(j)$. We call this the $j$th *significance graph*.

3. *Computing longest paths.* In each significance graph, we employ a depth-first search algorithm to explore all significant paths

$$\pi = (R_1, R_2, \ldots, R_m),$$

that is, sequences of vertices that are:

    (a) all significant, $s(R_k) = 1$;
    (b) all connected, $(R_k, R_{k+1}) \in \mathcal{E}(j)$.

We record the maximum path length in each significance graph:

$$L_{n,j}^{\max} = \max\{\text{length}(\pi) : \pi \text{ is a significant path in } \Sigma(j)\},$$

$$L_n^{\max} = \max_j L_{n,j}^{\max}.$$

4. *Decision.* We compare $L_n^{\max}$ with a length threshold:

    If $L_n^{\max} \leq L_n^*$, accept $\mathbf{H}_0$;    if $L_n^{\max} > L_n^*$, reject $\mathbf{H}_0$.

This defines the test, except for the specification of the thresholds $L_n^*$ and $N^*$. Asymptotic formulas for these thresholds will be given in Section 4.

A worked example of the multiscale significant-runs algorithm is illustrated in Figure 6. For this example the synthetic point cloud is mostly uniformly distributed, with fraction $\varepsilon \approx 1/20$ of points lying on a fixed smooth curve. The significance threshold was $N^* = 8$. For this choice of threshold, we have (under $\mathbf{H}_0$) $P\{N(R) > N^*\} \approx 0.00024$.

The longest run in this example has length 5, which in this case exceeds the run-length threshold $L_n^*$ and leads to rejection of the null hypothesis. For this small-sample setting, the threshold $L_n^* = 3$ was obtained by simulation rather than asymptotic theory. (Under the null hypothesis, we conducted 1200 experiments. In each experiment, a point cloud was generated and $L_n^{\max}$ was computed. The frequencies of $L_n^{\max} = 1, 2, 3, 4$ and 5 were 302, 873, 24, 0 and 1, respectively. Based on these results, $L_n^*$ can be set at either 2 or 3, giving a test with empirical level $P\{L_n^{\max} > 2\} = 25/1200$ or $P\{L_n^{\max} > 3\} = 1/1200$, resp.)

Asymptotic theory gives $L_n^* \approx 3.74$, which leads to the same decision.

For comparison, Figure 7 gives a simulated example in the null case $\varepsilon = 0$; the longest run has length 2 in this case. This simulation was typical; in the null example, rarely does the longest run exceed 2. Several properties of the algorithm are immediate:

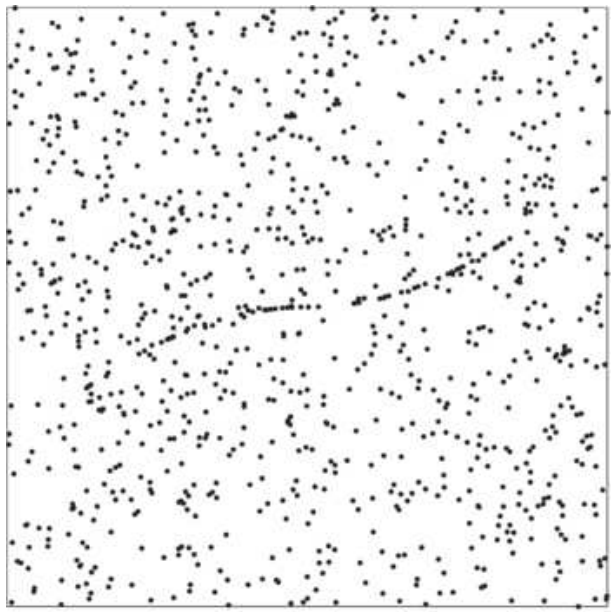

(a)

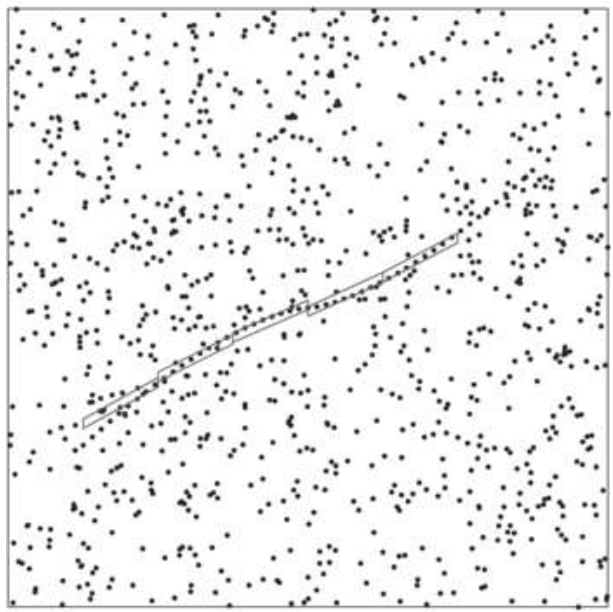

(b)

FIG. 6. *A uniform random scatter contaminated by $\varepsilon = 1/20$ points on a curve, together with the identified significant run (consisting of five strips); $n = 1024$.*

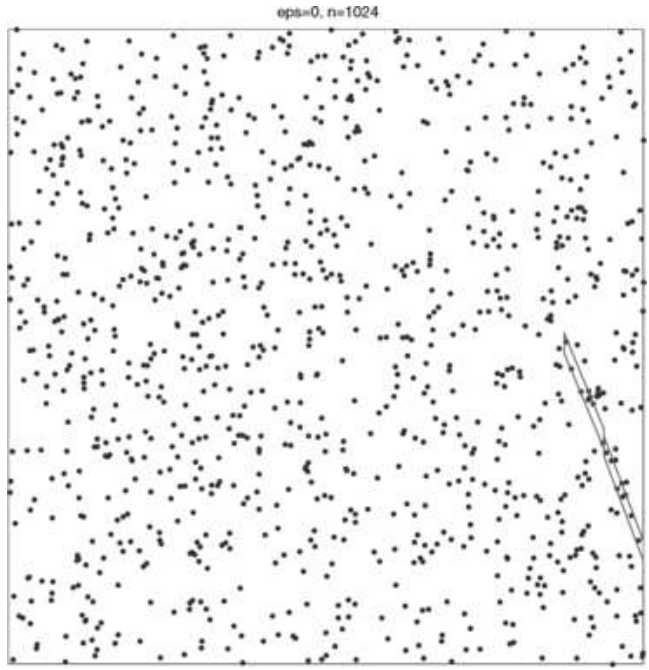

FIG. 7. *A uniform random scatter, together with the identified longest run of significant strips. At length 2 it is not a significant run; $n = 1024$.*

- *Complexity of strip counts.* The algorithm calculates all the $N(R)$ for all the anisotropic strips. This takes $O(n^2 \log(n))$ flops, where $n$ is the number of points sampled. Indeed, since each data point can belong to order $O(n \log(n))$ strips, by simply calculating which $R \ni X_i$ and incrementing a counter for those $R$, we get all $N(R)$.

- *Complexity of longest path.* The algorithm calculates the longest path in each significance graph. This takes work comparable to $O(n \log(n))$, based on depth-first search [2].

- *Storage requirement.* The algorithm stores all of the $N(R)$ counts and all the significance coefficients. This requires $O(n \log(n))$ storage.

To state our main result, we amend our notion of alternative hypothesis. For $\alpha > 1$, let Hölder$(\alpha, \beta, S)$ denote the collection of functions $f \in$ Hölder$(\alpha, \beta)$ with $|f'|_\infty \leq S$. Define

$$\mathbf{H}_1(\alpha, \beta, S, \tau) : X_i \overset{\text{i.i.d.}}{\sim} (1 - \varepsilon_n) \text{Uniform}(0,1)^2 + \varepsilon_n \text{Uniform}(\text{graph}(f)),$$

$$f \text{ in Hölder}(\alpha, \beta; S), \varepsilon_n > \tau \cdot n^{-\alpha/(1+\alpha)}.$$

THEOREM 3.1. *There is a single choice of thresholds $N^*$ and $(L_n^*)_n$ so that for every $\alpha \in (1, 2]$ and $\beta > 0$, there is $T_*(\alpha, \beta, S)$ with, for each $\tau > T_*$,*

$$P\{test\ rejects\ \mathbf{H}_0 | \mathbf{H}_1(\alpha, \beta, S, \tau)\} \to 1 \qquad as\ n \to \infty;$$

*at the same time*

$$P\{test\ rejects\ \mathbf{H}_0 | \mathbf{H}_0\} \to 0 \qquad as\ n \to \infty.$$

In words, under $\mathbf{H}_0$ the longest significant path is overwhelmingly unlikely to be substantially longer than $L_n^*$ for large $n$, while under each indicated $\mathbf{H}_1$ the longest significant path is overwhelmingly likely to be substantially longer than $L_n^*$. The threshold $T_*$ is within a constant factor $T_*/T_-$ of the optimal detection threshold; this shows that, up to constants, we can adaptively test for the existence of fragments of $C^\alpha$ graphs.

**4. Asymptotics of thresholds.** For practical finite-sample application of the test just described, it is of course possible to calibrate thresholds by conducting simulation experiments. For the proof of our main result, we specify thresholds in closed form. These thresholds are very conservative and our closed-form analysis yields vastly overstated estimates of error probabilities, which are nevertheless good enough for our proofs.

4.1. *Specification of $N^*$.* Fix a probability $p_0 < 1$. Define a counting threshold $N^+(\varepsilon, \lambda)$ with the property that

$$P\{\text{Poisson}(\lambda) > N^+\} \le \varepsilon,$$

where $\text{Poisson}(\lambda)$ denotes a Poisson random variable with parameter $\lambda$. Let $\text{Bin}(n, p)$ denote a binomial random variable with parameters $n$ and $p$. By Poisson approximation to the binomial, we have

$$P\{\text{Bin}(n, \lambda/n) > N^+(\varepsilon, \lambda)\} \le 2\varepsilon, \qquad n \ge n_0.$$

Set then $N^* = N^+(p_0/162, 2)$. Our definition of $N^*$ gives us the key property

(4.1)    $P\{s(R) = 1 | \mathbf{H}_0\} = P\{\text{Bin}(n, 2/n) > N^*\} \le p_0/81, \qquad n \ge n_0.$

4.2. *Specification of $L_n^*$.* We use a convenient, but nonstandard, notation borrowed from Arratia and Waterman [6]. For $0 < p < 1$, $\log_p$ denotes, in this unconventional notation, the *logarithm with base $1/p$*. At the same time, we maintain the original convention that if $b > 1$, $\log_b$ means log to the base $b$. By this convention, $\log_2$ is the traditional logarithm base two and $\log_{1/2}$ is actually the same quantity.

Define the length threshold

(4.2)    $$L_n^* = 3 \log_{p_0}(n),$$

where $p_0 \in (0, 1)$ is the same as in the specification of $N^*$. The underlying rationale for this choice is the Erdös–Rényi law (see [6]), which says:

> *In a sequence of $m$ i.i.d. Bernoulli random variables with probability $p$ of heads,*
> *the length of the longest run of pure heads $\sim \log_p(m)(1 + o_P(1))$.*

In effect, our definition makes $L_n^*$ very substantially longer than the length of the longest run of pure heads in a linear sequence of $O(n^{1/(1+\alpha)})$ coin tosses of a $p_0$ coin.

4.3. *Specification of $T_*$.* Associated to the parameters $N^*$ and $L_n^*$ will be the threshold $T_*$ at which $\mathbf{H}_1$ becomes detectable. To define that, set $p_1$ sufficiently close to 1 so that, for all $\alpha \in (1, 2]$ and for some $n_1 > 0$,

$$(4.3) \qquad \log_{p_1}(n^{1/(1+\alpha)}) \geq 2 \cdot L_n^*, \qquad n > n_1$$

($p_1 = p_0^{1/18}$ works). Implicitly, this choice again refers to the Erdös–Rényi law and, in some way, guarantees that under $\mathbf{H}_1$ there will be very long runs.

Define an intensity threshold $\Lambda^+(\varepsilon)$ with the property that

$$P\{\text{Poisson}(\Lambda^+(\varepsilon)) < N^*\} \leq \varepsilon.$$

Set $\lambda^* = \Lambda^+((1 - p_1)/2)$ and set

$$T_*(\alpha, \beta, S) = 2\lambda^* \cdot \beta^{1/(1+\alpha)} \cdot \sqrt{1 + S^2}.$$

We will show that for

$$(4.4) \qquad \varepsilon_n > T_* \cdot n^{-\alpha/(1+\alpha)}$$

the hypothesis $\mathbf{H}_1(\alpha, \beta)$ becomes detectable by our proposal. The central point will be the property

$$(4.5) \qquad n \cdot \varepsilon_n \cdot w(j^*(\alpha, \beta, n))/\sqrt{1 + S^2} \geq \lambda^*.$$

To prove this, we inspect the proof of Lemma 2.1 in the Appendix and note that, by definition of $T_*(\alpha, \beta, S)$ above and using notation from the Appendix,

$$n \cdot \varepsilon_n \cdot w(j^*(\alpha, \beta, n))/\sqrt{1 + S^2}$$
$$\geq n^{1/(1+\alpha)} \cdot T_* \cdot w(j^+(\alpha, \beta, n))/(2\sqrt{1 + S^2}) = \lambda^*.$$

Incidentally, we do not claim that the algorithm fails for $\tau < T_*$, only that we can prove it succeeds for $\tau > T_*$.

**5. Behavior under $\mathbf{H_1}$.** Let $f$ be the function in Hölder$(\alpha, \beta)$ that generates the curve that carries the fraction $\varepsilon_n$ of data. From Lemmas 2.1–2.3, we let $j = j^*$ and consider the tube $\mathcal{T}_j(f)$. For each region $R$ in this tube,

$$N(R) \sim \text{Bin}(n, (1 - \varepsilon_n)\text{area}(R) + \varepsilon_n \gamma(f, R)),$$

where $\text{area}(R) = 2/n$ and $\gamma$ denotes the relative arc length in the graph of $f$, obeying

$$\gamma(f, R) = \frac{\text{length}(\text{graph}(f|I))}{\text{length}(\text{graph}(f))} \geq \frac{w}{\sqrt{1 + S^2}};$$

here $I$ denotes the projection of $R$ on the $x$-axis.

By Poisson approximation to the binomial,

$$N(R) \overset{\text{approx.}}{\sim} \text{Poisson}(\mu),$$

where, using (4.5),

$$\mu \geq 1 + n\varepsilon_n w / \sqrt{1 + S^2}$$

$$\geq \lambda^*.$$

Hence, for every $R$ in this tube, we have, for all sufficiently large $n > n_3$,

(5.1) $$P\{N(R) > N^*\} \geq p_1.$$

Label the sequence of strips $R$ in this tube $R_0, \ldots, R_{w^{-1}-1}$. We want to know the probable length of the longest run of the form

$$N(R_k) > N^*, \ldots, N(R_{k+L}) > N^*.$$

Then, if $L_n$ is the length of the longest run in this special sequence, it follows that the longest run statistic we are computing over the entire graph can only be larger:

$$L_{n,j}^{\max} \geq L_n.$$

To show that the test rejects $\mathbf{H_0}$, we will show that

(5.2) $$P\{L_n > L_n^*\} \to 1, \qquad n \to \infty.$$

If we define

$$Z_i = \mathbb{1}_{\{N(R_k) > N^*\}},$$

we note that each $Z_i$ is Bernoulli with probability $p_i$, while

$$p_i \geq p_1.$$

We let $m = w^{-1} \geq \text{Const} \cdot n^{1/(1+\alpha)}$ and $p = p_1$, and we get, by the Erdös–Rényi law,

$$L_n > \log_{p_1}(n^{1/(1+\alpha)})(1 + o_P(1)).$$

However, by hypothesis we have chosen $p_1$ in such a way that

$$\log_{p_1}(n^{1/(1+\alpha)}) \geq 2L_n^*.$$

Hence (5.2) follows.

**6. Behavior under $\mathbf{H_0}$.** We need to show that, with overwhelming probability under $\mathbf{H}_0$, there will be no runs in the graph that exceed $L_n^*$. We start by arguing that

$$P\{\text{a significant path of length } L \text{ starts at given } R|\mathbf{H}_0\} \leq p_0^L.$$

Indeed, by choice of $N^*$, for each $R$,

$$P\{s(R) = 1|\mathbf{H}_0\} = P\{\text{Bin}(n, 2/n) > N^*\} \leq p_0/81.$$

Now each region $R$ has 81 neighbors in $\mathcal{G}(j)$, and so by Boole's inequality,

$$P\{s(R') = 1 \text{ for at least one neighbor of } R|\mathbf{H}_0\}$$
$$\leq \sum_{R' \in \text{neighbors}(R)} P\{s(R') = 1|\mathbf{H}_0\}$$
$$\leq \# \text{neighbors} \cdot (p_0/81)$$
$$= 81 \cdot (p_0/81) = p_0.$$

Now if we are looking for a significant path of length greater than $L$, we need that starting at some vertex, it has $s(R) = 1$ and is connected to a vertex with $s(R') = 1$, et cetera. For a given starting point $R$, the probability of this event is bounded by $p_0^L$ via negative correlation.

We now note there are at most $M_j = w^{-1} \cdot \delta_1^{-1} \cdot \delta_2^{-1} \cdot 2S$ starting points for paths in $\mathcal{G}(j)$. By Boole's inequality,

$$P\{\text{there is a significant path of length } L \text{ occurring in } \mathcal{G}(j)|\mathbf{H}_0\}$$
$$\leq \# \text{ (starting points at level } j)$$
$$\times P\{\text{significant path of length } L \text{ starting at } R|\mathbf{H}_0\}$$
$$\leq M_j \cdot p_0^L.$$

Take $\log_2$:

$$\log_2(M_j) + \log_2(p_0) \cdot L = \log_2(w^{-1} \cdot \delta_1^{-1} \cdot \delta_2^{-1} \cdot 2S) + \log_2(p_0) \cdot L$$
$$= 2J - 2j + 3 + \log_2(S) + \log_2(p_0) \cdot L$$
$$\leq 2J + C + \log_2(p_0) \cdot L.$$

For $L = L_n^*$, the last expression on the right-hand side tends to $-\infty$ as $n$ increases. Hence

$$P\{\text{there is a run of length } L_n^*|\mathbf{H}_0\} \to 0, \qquad n \to \infty.$$

**7. Discussion.** We have given only a sampling of results in a specific problem of geometric detection; much more could be done. At the same time, our results are closely related to many ideas in the literature of computer vision. We briefly indicate possible variations and sketch a few such connections.

7.1. *Variations on the filament model.* We have only discussed a subset of what could pass for filamentary structure in point-cloud data. Other notions of filamentarity include (a) curves that have less regularity than one derivative, (b) curves that have more regularity than two derivatives and (c) curves that are not describable as simple graphs $(x, y = f(x))$. Our aim in this paper is to stimulate discussion; we believe that all such generalizations will be of interest in appropriate applications areas. We make a few brief remarks.

- *Curves that have regularity $\alpha \leq 1$.* If we consider curves $(x, f(x))$, where $f$ is Hölder-$\alpha$ with $\alpha \leq 1$, we are considering curves without tangents. We thus discard completely the notion of good continuation based on alignment of tangents. In designing graph-based detectors, it is only required to use axis-oriented rectangular regions, so that only position (not slope) matters, and in which the rectangles are now taller than they are wide; Connectivity involves only position (not orientation). The statistical treatment based on graphs and runs turns out to be the same; the graphs simply have less structure because proximity does not involve similarity of slopes.

- *Curves that have regularity $\alpha > 2$.* It makes sense to ask about smoothness of higher order, for example, to consider $2 < \alpha \leq 3$. By [5], there will continue to be about $n^{1/(1+\alpha)}$ points on some Hölder curve, even for $\alpha$ in this range. To sensitively detect such higher-order smoothness would require a different set of regions than the one discussed here—curved ones with parabolic midlines—and a notion of good continuation based on matching of sides, and matching of slopes and curvatures of midlines. Preliminary calculations suggest that analogous "adaptivity" results hold in such a setting. Related discussion can be found in [3, 5].

- *Curves that are not graphs.* The Introduction suggested that the approach described here can be adapted to detect general plane curves, that is, curves that are not graphs. This adapts ideas from our work on beamlet graphs [4, 16]. We define a family of directed graphs based on regions modeled on "dyadically thickened beamlets" with various degrees of thickening. In this directed graph structure, strips can have all orientations, including vertical and horizontal, so that the graph constraint is removed. Connectivity between beamlets is based once more on good continuation principles—in this case, continuation of plane polygons rather than polygonal graphs. Otherwise the algorithms are identical. More details can be found in [5], where this structure is utilized to prove a theoretical result.

Still other variations on our model are possible. The idea that data are uniformly sampled from a curve of zero width can be varied in several ways:

- *Nonuniform sampling along filaments.* A referee suggested that rather than from a uniform density, data might arise instead from a density that is bounded away from zero and infinity. The same methodology developed here works in that case without change, except that the analysis under $\mathbf{H}_1$ seemingly becomes more involved.

- *Finite thickness.* A referee also suggested that rather than from a curve of zero thickness, the data might arise instead from a tube of finite thickness. The methodology developed here works without change in such a case, but the model $\mathbf{H}_1$ is different and the statement of results becomes different. Thus, if the thickness of the tube is finite, but smaller than the width of the regions that are adapted to the underlying filament, the detectability results are the same as here. If the width is greater than that of regions adapted to the curve, then the detection threshold becomes higher and it takes systematically larger numbers of points near a curve to reject $\mathbf{H}_0$.

- *Finite resolution.* A referee also suggested that the data might be of finite accuracy, for example, either subject to rounding or to noise. In either case, the situation is much the same as in the immediately preceding comment. If the inaccuracy is small compared to the width of the optimally fitted regions, its impact is negligible. If the inaccuracy is larger than that width, then the detection threshold becomes higher and rejection of $\mathbf{H}_0$ requires systematically larger numbers of points near a curve.

We leave exploration of such issues for future work.

7.2. *Beyond detection.* A referee pointed out the desirability of not just detecting the presence of filamentary structures, but also estimating the detailed location and shape of any filaments that are detected. Of course, conceptually detection and estimation are quite different tasks; moreover, unless detection is possible, estimation is impossible. Empirically, our methodology actually provides an estimator when the filament is not just detected, but strongly detected. As the number of points sampled from the filament increases well beyond the detection threshold, simulations show that the longest run becomes overwhelmingly likely to trace out a series of regions that bracket the curve tightly. It seems likely that this effect could be proven rigorously. However, it seems rather delicate to formulate and study an appropriate notion of asymptotically efficient estimation.

7.3. *Extensions beyond the filament model.* We can extend beyond the setting of two-dimensional point clouds in at least three ways: going to higher dimensions, observing vectors rather than points and observing pixel imagery on a grid rather than scattered points.

7.3.1. *Structures in higher dimensions.* The analogous detection problem in $d$-dimensional space—finding a curve or surface that contains an unexpectedly large number of points—has been considered by the authors in [3, 5]. That work provides nonadaptive detectors, that is, detectors that assume knowledge of the Hölder class.

The ideas developed in this paper can be directly applied to multiscale detection of filamentary structure in $d$-dimensional point clouds, $d > 2$. A more ambitious generalization—detection of codimension-$k$ surfaces in $d$-dimensional point clouds—seems possible, but also messier. For example, for $d = 3$ and $k = 1$, we are attempting to find a surface in $d$-dimensional space that contains an inordinately large number of points. The nodes of each anisotropy graph are planar slabs of volume $2/n$ and the neighborhood structure in the graph is, while conceptually analogous to the case considered here, far more complex to write out. The comparable adaptive detection theorems hold in that setting, although we omit details.

7.3.2. *Vector fields.* Suppose that instead of data on points $X_i$, we have data on tangent vectors, that is, on pairs $(X_i, \theta_i)$ that name both a position and a direction. Experiments in perceptual psychophysics (e.g., [19, 28]) suggest that this is a much more potent stimulus to "curve finding" than simply the display of random dots. Biological evidence about early vision suggests that individual receptive fields fire when both location and orientation offer matches.

As a null hypothesis, we suppose that the $X_i$'s are i.i.d. uniform$[0, 1]^2$, while the $\theta_i$'s are i.i.d. uniform$[-\pi/2, \pi/2]$. As an alternative hypothesis, we could posit that a small fraction of the $X_i$ lie on a curve, $X_i = (x_i, f(x_i))$, and the $\theta_i$ specify angles parallel to the line with slope $f'(x_i)$.

Our paper [5] showed that such tangent vector data are substantially more powerful for identifying filamentary structure than the point data discussed so far in this paper. Namely, such data give us the ability to detect filaments that contain much smaller fractions of data points. In fact, a reliable test for a Hölder$(2, 1)$ filament can be based on agreement with more than $Tn^{1/4}$ tangent vectors, rather than $Tn^{1/3}$ points.

The multiscale data structures used in this paper can be applied in the tangent vector setting, where we match $(X_i, \theta_i)$ to a region based not only on $X_i \in R'$, but also on $\theta_i$ matching the slope of the midline of $R$. The resultant algorithm can take advantage of this more stringent matching structure to speed up the counting process, because each tangent vector will lie in only one region in a given anisotropy class, resulting in $O(\log(n))$ flops per data point rather than $O(n \log(n))$. Searching for significant paths will be significantly faster as well, since there can only be $n \log(n)$ starting places for a path. Hence the whole algorithm can run in $O(n \log(n))$ flops. An analysis that parallels the one given here shows that a multiscale multianisotropic significance runs algorithm can provide a detection threshold that is optimal to within a constant factor. Huo and co-workers [23] gave more details from the computational aspect of this problem.

7.3.3. *Pixel imagery.* In another direction, we might consider data types used to model digital imagery, for example, arrays $(y(i_1, i_2) : 0 \leq i_1, i_2 < n)$, where

$$y_{i_1, i_2} = \xi_f(i_1, i_2) + \sigma z(i_1, i_2), \qquad 0 \leq i_1, i_2 < n,$$

$(z(i_1, i_2))$ is a Gaussian white noise, $\sigma$ is the noise level and $\xi_f$ is a pixel array with nonzero values only on pixels that intersect graph($f$). Despite appearances, this problem is closely related to the present problem and analogous detection theorems are true.

In effect, we form a family of anisotropic multiscale strips and sum the pixels that intersect those strips, producing detector statistics $X(R)$ that can be significance-tested in a way that parallels the counts $N(R)$ considered in this paper, only with Gaussian rather than Poisson threshold analysis. The underlying data structures and arguments can be understood as, roughly speaking, a mixture of the ideas of this paper and those of another paper by these authors [4]. In particular, the strips considered here were called axoids in that work. More details can be found in [21, 22]. In the digital array setting, there is a fast beamlet transform to rapidly compute all the required $X(R)$ detector statistics.

7.4. *On the uniform background assumption.* A referee remarked that many practical imaging problems do not involve small departures from a uniform background. This is no doubt true. However, there is an ever-expanding array of imaging problems, and some look for changes between one image and another, for example, in studying arterial blood flow or in change detection in scene surveillance. In the case of no change, the background is quite literally uniform. Also, many problems with nonuniform background are transformable to uniform background; thus edge detectors and object detectors are typically operated at so-called constant false alarm rate. Then they give, under the null hypothesis of no object present, roughly a constant number of events per unit area, and so the constant false alarm rate transformation forces a uniform background under $\mathbf{H}_0$. Finally, the intellectually important issues explored here seem clearest in the uniform case and the data structures developed here are known to be useful in more complex settings [4].

The main point, however, is well taken—detection of objects against nonuniform clutter rather than uniform scatter remains a challenging area for further work.

7.5. *Relationships to other work.*

7.5.1. *Parametric detection.* In this paper we considered detection of points on nonparametric curves. In certain cases, one is interested in points along lines [8] or on parabolas [1]. For an attractive nonmultiscale approach to such detection problems, see [11, 12, 13] and [33], Section 6.3.

In the authors' paper [4], which was just mentioned, it was shown how to develop multiscale geometric detectors for line segments in digital imagery, for parametric forms such as circles, rectangles and ellipses. Those ideas could be adapted to the present setting to find situations where data have an elevated density over a blob or along some line. In the end, the underlying computations involve dyadic multiscale rectangles and strips, and the ideas are closely related to those in this paper.

7.5.2. *Multiscale geometric analysis.* The tools described here are closely related to a variety of tools in multiscale methods; see [10, 15, 16, 18, 24, 27, 32] for discussions of related tools applied in image analysis and in mathematical analysis. We differ here in our use of a multiscale multianisotropy collection of analyzing regions that is organized and exploited in a specific way and for a specific purpose.

7.5.3. *Object grouping and neural architecture.* In the literature of computer vision, there is extensive discussion of object grouping in perception [7]. The problem considered here—recognizing a curve against a background of random points—fits in this tradition. There are even experiments in psychophysics that test the ability of the human visual system to accomplish similar tasks [19, 25]. In effect, what we have discussed here—(near-) optimal detection—corresponds to what psychophysicists call "the ideal observer" [26]. In this connection, we have exhibited a simple multiscale architecture that can provide a near-optimal detector for a very wide range of stimuli—curves of any of a wide range of degrees of smoothness.

In the Gestalt theory of perception, there arises the concept of good continuation [19, 35]; experiments show that the visual system will respond better to curvilinear stimuli that follow a good continuation of an initial pattern [30]. Here we have operationalized this principle by a regular connectivity pattern in a specific graphical structure. We have shown that with this implementation, we get a near-optimal detector, thus validating the significance of such a good continuation principle. Note well that the connectivity pattern is invariant; that is, it applies the same way at all nodes of the graph.

This architectural simplicity is striking when compared with the vast speculative literature that proposes "neural architectures" for visual perception. What we have shown is that by starting from a large collection of elements that are sensitive to (i.e., accumulate counts in) receptive fields at a variety of lengths, widths, orientations and locations, and then connecting such elements to other elements by a simple invariant rule, one very sensitively recognizes the existence of curvilinear stimuli simply by the existence of long connected paths.

Perhaps this bears comparison with biological evidence. Functional magnetic resonance imaging studies [9] suggest that there are centers in primate brains that seem responsible for integrating local information into recognition of long curvilinear structures [29]. It would be interesting to know whether such integration has any resemblance to the simple multiscale connection mechanism employed here.

## APPENDIX

PROOF OF LEMMA 2.1.   Extend notation so that $w(j) = 2^{-j}$ can have both real and integer arguments, and $t(j) = 2^{-(J-j)+1}$ as well. Let $j^+ = j^+(\alpha, \beta, n)$ satisfy

$$2\beta w(j^+)^\alpha = t(j^+),$$

that is,

$$2\beta 2^{-\alpha j^+} = 2^{-(J-j^+)+1}.$$

Let $j^* = \lceil j^+ \rceil$ be the next larger integer, so that $w(j^+)/2 \leq w(j^*) \leq w(j^+)$ and $t(j^+) \leq t(j^*) \leq 2t(j^+)$. Then because $1 \leq \alpha \leq 2$, $2^\alpha \leq 4$ and so

$$2\beta w(j^*)^\alpha \leq 2\beta w(j^+)^\alpha = t(j^+) \leq t(j^*)$$
$$\leq 2 \cdot t(j^+) = 4\beta w(j^+)^\alpha \leq 16\beta w(j^*)^\alpha.$$

Now let $w = w(j^*)$ and $t = t(j^*)$, and substitute in the last display, getting (2.1).   □

PROOF OF LEMMA 2.2.   Remember that $f \in \text{Hölder}(\alpha, \beta)$ satisfies $f : [0, 1] \to [0, 1]$ and

$$|f'(x) - f'(y)| \leq \alpha\beta|x - y|^{\alpha-1}.$$

We saw that this implies

(A.1)        $|f(x) - f(y) - f'(y)(x - y)| \leq \beta|x - y|^\alpha, \qquad x, y \in [0, 1].$

To prove (2.2), we use the notation $\dot{f}_k(x)$ for the affine function tangent to $\text{graph}(f)$ at $x_k$. Figure 4 illustrates our notation.

Using (A.1), with $I_k$ denoting $[kw, (k+1)w)$,

$$|f(x) - \dot{f}_k(x)| \leq \beta(w/2)^\alpha \leq t/4, \qquad x \in I_k.$$

We also note that if $\dot{g}_k(x) = l_1\delta_1 + l_2\delta_2(x - x_k)$, then

$$|\dot{f}_k(x) - \dot{g}_k(x)| \leq |f(x_k) - l_1\delta_1| + |f'(x_k) - l_2\delta_2||x - x_k|$$
$$\leq \delta_1/2 + \delta_2/2 \times w/2$$
$$\leq t/8 + t/16.$$

We conclude that

$$(A.2) \qquad |f(x) - \dot{g}_k(x)| \leq t/2, \qquad x \in I_k.$$

On the other hand, the region $R(j, k, \ell_1, \ell_2)$ has $\dot{g}_k(x)$ as its midline and is of half-height $t/2$. The desired relationship (2.2) follows. $\square$

PROOF OF LEMMA 2.3. We use the same notation as in the proof of Lemma 2.2, as illustrated in Figure 4.

It is enough to show that

$$(A.3) \qquad |\dot{g}_{k+1}(x_{k+1}) - \dot{g}_k(x_{k+1})| \leq t$$

and

$$(A.4) \qquad |\dot{g}'_{k+1}(x_{k+1}) - \dot{g}'_k(x_{k+1})| \leq t/w.$$

It then follows that there is an edge in $\mathcal{E}(j)$ that connects $R(j, k, \ell_1, \ell_2)$ to $R(j, k+1, \ell'_1, \ell'_2)$, where $\ell'_1$ and $\ell'_2$ are the values associated to $\dot{g}_{k+1}$.

The following inequalities flow either from Hölder conditions or from simple rounding involved in quantization:

$$|\dot{g}_{k+1}(x_{k+1}) - f(x_{k+1})| \leq \delta_1/2 = t/8,$$
$$|f(x_{k+1}) - \dot{f}_k(x_{k+1})| \leq \beta w^\alpha \leq t/2,$$
$$|\dot{f}_k(x_{k+1}) - \dot{g}_k(x_{k+1})| \leq \delta_1/2 + \delta_2/2 \cdot w = t/4.$$

Combining these with the triangle inequality yields (A.3). Similarly,

$$|\dot{g}'_{k+1}(x_{k+1}) - f'(x_{k+1})| \leq \delta_2/2 = t/(8w),$$
$$|f'(x_{k+1}) - \dot{f}'_k(x_{k+1})| \leq \alpha\beta w^{\alpha-1} \leq t/(2w),$$
$$|\dot{f}'_k(x_{k+1}) - \dot{g}'_k(x_{k+1})| \leq \delta_2/2 = t/(8w);$$

combining these identities gives (A.4). $\square$

**Acknowledgments.** We would like to thank Emmanuel Candès, Hagit Hel-Or, Aapo Hyvärinen, Jean-Luc Starck and Brian Wandell for helpful discussions and references.

## REFERENCES

[1] ABRAMOWICZ, H., HORN, D., NAFTALY, U. and SAHAR-PIKIELNY, C. (1997). An orientation-selective neural network for pattern identification in particle detectors. In *Advances in Neural Information Processing Systems* **9** (M. Mozer, M. I. Jordan and T. Petsche, eds.) 925–931. MIT Press, Cambridge, MA.

[2] AHO, A. V., HOPCROFT, J. E. and ULLMAN, J. D. (1983). *Data Structures and Algorithms.* Addison–Wesley, Reading, MA. MR0666695

[3] ARIAS-CASTRO, E. (2004). Graphical structures for geometric detection. Ph.D. dissertation, Stanford Univ.

[4] ARIAS-CASTRO, E., DONOHO, D. L. and HUO, X. (2005). Near-optimal detection of geometric objects by fast multiscale methods. *IEEE Trans. Inform. Theory* **51** 2402–2425.

[5] ARIAS-CASTRO, E., DONOHO, D. L., HUO, X. and TOVEY, C. (2005). Connect-the-dots: How many random points can a regular curve pass through? *Adv. in Appl. Probab.* **37** 571–603. MR2156550

[6] ARRATIA, R. and WATERMAN, M. S. (1989). The Erdös–Rényi strong law for pattern matching with a given proportion of mismatches. *Ann. Probab.* **17** 1152–1169. MR1009450

[7] BUHMANN, J. M., MALIK, J. and PERONA, P. (1999). Image recognition: Visual grouping, recognition, and learning. *Proc. Natl. Acad. Sci. USA* **96** 14,203–14,204.

[8] COPELAND, A. C., RAVICHANDRAN, G. and TRIVEDI, M. M. (1995). Localized Radon transform-based detection of ship wakes in SAR images. *IEEE Trans. Geoscience and Remote Sensing* **33** 35–45.

[9] COURTNEY, S. M. and UNGERLEIDER, L. G. (1997). What fMRI has taught us about human vision. *Current Opinion in Neurobiology* **7** 554–561.

[10] DAVID, G. and SEMMES, S. (1993). *Analysis of and on Uniformly Rectifiable Sets.* Amer. Math. Soc., Providence, RI. MR1251061

[11] DESOLNEUX, A., MOISAN, L. and MOREL, J.-M. (2000). Meaningful alignments. *Internat. J. Computer Vision* **40** 7–23.

[12] DESOLNEUX, A., MOISAN, L. and MOREL, J.-M. (2003). A grouping principle and four applications. *IEEE Trans. Pattern Analysis and Machine Intelligence* **25** 508–513.

[13] DESOLNEUX, A., MOISAN, L. and MOREL, J.-M. (2003). Maximal meaningful events and applications to image analysis. *Ann. Statist.* **31** 1822–1851. MR2036391

[14] DONOHO, D. L. (1997). CART and best-ortho-basis: A connection. *Ann. Statist.* **25** 1870–1911. MR1474073

[15] DONOHO, D. L. (1999). Wedgelets: Nearly minimax estimation of edges. *Ann. Statist.* **27** 859–897. MR1724034

[16] DONOHO, D. L. and HUO, X. (2002). Beamlets and multiscale image analysis. In *Multiscale and Multiresolution Methods. Lecture Notes Comput. Sci. Eng.* **20** 149–196. Springer, Berlin. MR1928566

[17] DONOHO, D. L. and JOHNSTONE, I. M. (1995). Adapting to unknown smoothness via wavelet shrinkage. *J. Amer. Statist. Assoc.* **90** 1200–1224. MR1379464

[18] DONOHO, D. L. and LEVI, O. (2004). Fast X-ray and beamlet transforms for three-dimensional data. In *Modern Signal Processing* (D. N. Rockmore and D. M. Healy, Jr., eds.) 79–116. Cambridge Univ. Press. MR2075950

[19] FIELD, D., HAYES, A. and HESS, R. (1993). Contour integration by the human visual system: Evidence for a local "association field." *Vision Research* **33** 173–193.

[20] Ho, M.-W. (2004). In search of the sublime. Institute of Science in Society. Available at www.i-sis.org.uk/sublime.php.

[21] Huo, X., Chen, J. and Donoho, D. L. (2003). Multiscale detection of filamentary features in image data. In *Wavelets: Applications in Signal and Image Processing X* (M. A. Unser, A. Aldroubi and A. F. Laine, eds.) 592–606. SPIE, Bellingham, WA.

[22] Huo, X., Chen, J. and Donoho, D. L. (2003). Multiscale significance run: Realizing the "most powerful" detection in noisy images. In *Proc. Thirty Seventh Asilomar Conference on Signals, Systems, and Computers* **1** 321–326. IEEE, Piscataway, NJ.

[23] Huo, X., Donoho, D. L., Tovey, C. and Arias-Castro, E. (2004). Dynamic programming methods for "connecting the dots" in scattered point sets. Technical report, Dept. Statistics, Stanford Univ.

[24] Jones, P. W. (1990). Rectifiable sets and the traveling salesman problem. *Invent. Math.* **102** 1–15. MR1069238

[25] Kovacs, I. and Julesz, B. (1993). A closed curve is much more than an incomplete one: Effect of closure in figure-ground segmentation. *Proc. Natl. Acad. Sci. USA* **90** 7495–7497.

[26] Legge, G. E., Kersten, D. and Burgess, A. E. (1987). Contrast discrimination in noise. *J. Opt. Soc. Amer. A* **4** 391–404.

[27] Lerman, G. (2003). Quantifying curvelike structures of measures by using $L_2$ Jones quantities. *Comm. Pure Appl. Math.* **56** 1294–1365. MR1980856

[28] Levi, D. M. and Klein, S. A. (2000). Seeing circles: What limits shape perception? *Vision Research* **40** 2329–2339.

[29] Mendola, J. D., Dale, A. M., Fischl, B., Liu, A. K. and Tootell, R. B. H. (1999). The representation of illusory and real contours in human cortical visual areas revealed by functional magnetic resonance imaging. *J. Neuroscience* **19** 8560–8572.

[30] Pizlo, Z., Salach-Golyska, M. and Rosenfeld, A. (1997). Curve detection in a noisy image. *Vision Research* **37** 1217–1241.

[31] Qaddoumi, N., Ranu, E., McColskey, J. D., Mirshahi, R. and Zoughi, R. (2000). Microwave detection of stress-induced fatigue cracks in steel and potential for crack opening determination. *Research in Nondestructive Evaluation* **12** 87–103.

[32] Sharon, E., Brandt, A. and Basri, R. (2000). Fast multiscale image segmentation. In *Proc. IEEE Conference on Computer Vision and Pattern Recognition* **1** 70–77.

[33] Small, C. G. (1996). *The Statistical Theory of Shape.* Springer, Berlin. MR1418639

[34] Tupin, F., Maitre, H., Mangin, J.-F., Nicolas, J.-M. and Pechersky, E. (1998). Detection of linear features in SAR images: Application to road network extraction. *IEEE Trans. Geoscience and Remote Sensing* **36** 434–453.

[35] Wertheimer, M. (1938). *Laws of Organization in Perceptual Forms.* Harcourt Brace, London.

E. Arias-Castro
Department of Mathematics
University of California, San Diego
9500 Gilman Drive
La Jolla, California 92093-0112
USA
E-mail: eariasca@math.ucsd.edu

D. L. Donoho
Department of Statistics
Stanford University
Stanford, California 94305-4065
USA
E-mail: donoho@stat.stanford.edu

X. Huo
School of Industrial
  and Systems Engineering
Georgia Institute of Technology
Atlanta, Georgia 30332-0205
USA
E-mail: xiaoming@isye.gatech.edu